\documentclass[11pt]{article}

\usepackage{graphicx}
\usepackage{amsthm}
\usepackage{amsmath}
\usepackage{amssymb}
\usepackage{epstopdf}

\theoremstyle{plain}
\newtheorem{theorem}{Theorem}
\newtheorem{lemma}[theorem]{Lemma}

\theoremstyle{definition}

\def\@seccntformat#1{Section\ \csname the#1\endcsname\quad}

\def\be{\begin{enumerate}}
\def\ee{\end{enumerate}}
\def\bi{\begin{itemize}}
\def\ei{\end{itemize}}

\title{A note on transversal knots which are closed 3-braids}
\author{Joan S. Birman and William W. Menasco}
\date{April 7, 2007}

\begin{document}
\maketitle

\begin{abstract}
The topological classification of knots that are closed 3-braids is shown to lead to a classification theorem for tranversal knots that are represented by closed 3-braids.  A list is given of all low crossing examples of transversally non-simple knots that are closed 3-braids.
\end{abstract}

A transversal knot type is said to be {\it transversally simple} if it is determined, up to transversal isotopy, by its topological knot type and its Bennequin or self-linking invariant.  In the absence of other invariants,  it was natural to ask was whether all transversal knot types are transversally simple.   In the manuscript  \cite{BM-stab-II} the authors of this note proved that the answer is no, by exhibiting an infinite family of pairs of transversal knot types, all closed 3-braids,  which are smoothly isotopic, but not transversally isotopic. Soon after  \cite{BM-stab-II} was posted, additional examples of the same type were found by Etnyre and Honda \cite{EH}.  The proofs in both \cite{BM-stab-II} and \cite{EH} were indirect, and did not lead to computable invariants. 

In the years since \cite{BM-stab-II} was posted there was new work, aimed at finding computable invariants of transversal knot type. 
In \cite{Plam1} Olga Plamenevskaya studied  knots via their 2-fold branched covers, hoping that techniques from Heegaard Floer Homology on the covering spaces would detect the differences between the knots.  In a second paper \cite{Plam2} Plamenevskaya developed a different, and now explicit invariant of transversal knot types, this time using Khovanov homology.  Her invariant was generalized by H. Wu in \cite{Wu}. However, a private communication with Professor Wu indicated that even the improved invariant failed to detect the examples in \cite{BM-stab-II}.  Most recently, an invariant of transverse knots was discovered by Ozsvath, Szabo and  Thurston in \cite{OST}.
In \cite{NOT}  Ng, Ozsvath and Thurston used the invariants in \cite{OST} to exhibit many pairs of transverse knots which are not transversally simple.   The question then arose as to whether the Ozsvath-Szabo-Thurston invariants, together with the classical invariants, sufficed to classify transversal knot types.   The answer seems to be no, because the techniques  in \cite{NOT} fail to distinguish at least one pair among the low crossing number examples in \cite{BM-stab-II}.  The situation therefore continues to be interesting. 

This flurry of activity, and in particular a lively interchange between the first author and Dylan Thurston when the work in \cite{NOT} was in progress, led us to go back and take another look at the work we had done in \cite{BM-stab-II}.  The purpose of this note is to restate the classification theorem for links that are closed 3-braids in a form that is, perhaps, more useful than that in \cite{BM-3-braids}, because it gives a complete enumeration of all the transversally non-simple knots which are represented by closed 3-braids, and a list of all low-crossing examples, whereas  in \cite{BM-stab-II} we only gave an infinite sequence of examples.  The complete enumeration is possible because of the results in  \cite{KL}. When \cite{KL} was written, the main result in \cite{BM-stab-I} but not \cite{BM-stab-II} was known.  

We remark that the main theorems in the papers \cite{BM-3-braids},\cite{BM-stab-I} and \cite{BM-stab-II} all had long and difficult proofs.  Moreover, the very useful paper \cite{KL} seems not to be widely-known.  Over the years, many people have asked us for a concise, explicit statement of  Theorem~\ref{T;classification of links that are closed 3-braids} below.  When invariants of transversal links first became an active subject, many people were interested in Theorem~\ref{T:transversal closed 3-braids}, however even the experts seemed to be unaware that there were low crossing number examples.  Therefore, even though the results in this note are very close to known results,  we hope that this note will make them more accessible.

We begin our work with a definition. 

{\bf Definition:} A closed 3-braid is said to {\it admit a flype} if it has one of the special forms sketched in Figure~1.  The associated closed braids in each pair are $\lambda = \sigma_1^u\sigma_2^v\sigma_1^w\sigma_2^\epsilon$ and 
$\lambda' = \sigma_1^u\sigma_2^\epsilon\sigma_1^w\sigma_2^v$.  The flype is  negative or positive, according as $\epsilon$ (the unique crossing in each picture) is $\mp$.  
\begin{figure}[htpb!]
\label{F:negative-flype}
\centerline{\includegraphics[scale=1.0] {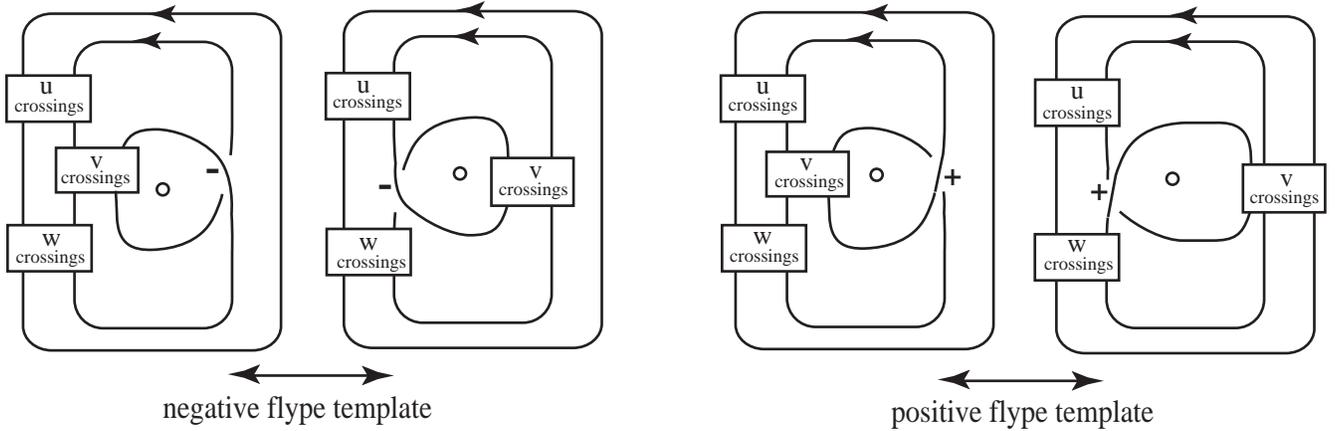}}
\caption{Negative and Positive 3-braid flypes}
\end{figure}

{\bf Example 1:}  The braids  $\sigma_1^3\sigma_2^4\sigma_1^5\sigma_2^{-1}$ and $\sigma_1^3\sigma_2^6\sigma_1^3\sigma_2^{-1}$ are in the same conjugacy class, and both admit negative flypes.  This is an example of a more general phenomenon, which had been overlooked in earlier work.  Let $\approx$ denote `is conjugate to':

\begin{lemma}
\label{L:conjugate triples}
If a 3-braid is in $\epsilon$-flype position $\sigma_1^u\sigma_2^v\sigma_1^w\sigma_2^\epsilon$, then so are the following other 3-braids, all of which are in the same conjugacy class and in $\epsilon$-flype position:  
\begin{equation}
\label{E:braid representation}
\sigma_1^u\sigma_2^v\sigma_1^w\sigma_2^\epsilon \approx \sigma_1^{v+\epsilon}\sigma_2^{w-\epsilon}\sigma_1^u\sigma_2^{\epsilon} \approx \sigma_1^w\sigma_2^{u-\epsilon}\sigma_1^{v+\epsilon}\sigma_2^\epsilon
 \end{equation}	
 \end{lemma}
 
{\bf Proof:} The relation $\sigma_i\sigma_{i+1}\sigma_i = \sigma_{i+1}\sigma_i\sigma_{i+1}$  in the braid group $B_n$ implies:
\begin{equation}
\sigma_i^\epsilon \sigma_{i+1}^k\sigma_i^{-\epsilon} = \sigma_{i+1}^{-\epsilon}\sigma_i^k\sigma_{i+1}^\epsilon
\end{equation}
Now consider the braid $\lambda = \sigma_1^u\sigma_2^v\sigma_1^w\sigma_2^\epsilon$.  We may rewrite $\lambda$ in another way that also shows it admits a negative flype:
$$ (\sigma_1^w)(\sigma_2^\epsilon\sigma_1^{u}\sigma_2^{-\epsilon})(\sigma_2^{v+\epsilon})= 
	(\sigma_1^w)(\sigma_1^{-\epsilon}\sigma_2^u\sigma_1^{\epsilon})(\sigma_2^{v+\epsilon})  \approx
	\sigma_2^{w-\epsilon}\sigma_1^u\sigma_2^{\epsilon}\sigma_1^{v+\epsilon}  \approx
	\sigma_1^{v+\epsilon}\sigma_2^{w-\epsilon}\sigma_1^u\sigma_2^{\epsilon}.$$
Applying this process again, the lemma is proved. Applying it a third time, we recover the braid $\lambda$. $\|$ 

 {\bf Remark:} We conjecture, but could not prove, that (\ref{E:braid representation}) gives all of the braid representatives in each conjugacy class which are in flype position and have minimal letter length.

\begin{theorem} {\bf The classification theorem for links that are closed 3-braids} {\rm (\cite{BM-3-braids}, \cite{KL})}
\label{T;classification of links that are closed 3-braids}
Let $\mathcal{L}$ be a link type which is represented by a closed 3-braid $\lambda$.  Then one of the following holds:
\be
\item [{\rm (1)}] $\lambda$ is conjugate to $\sigma_1^k\sigma_2^{\pm 1}$ for some $k\in\mathbb{N}$. In this case the closed braid associated to $\lambda$ represents either the unknot or a 2-component unlink or a torus link of type $(2,k)$, all of which have braid index less than 3.
\item  [{\rm (2)}]  $\lambda$ has braid index 3 and is represented by a unique conjugacy class of closed 3-braids.
\item [{\rm (3)}]  $\lambda$ has braid index 3, and is represented by exactly two distinct conjugacy classes of closed 3-braids. To enumerate the two conjugacy classes, when there are two, let $(u,v,w,\epsilon)$ be integers, with $\epsilon = \pm 1$.   Assume that $(u, v+\epsilon, w)$ are all distinct, also that $|v| \geq 2$ and that neither $u$ nor $w$ is equal to $0, \epsilon,  2\epsilon$.   Then the two conjugacy classes are represented by $\lambda = \sigma_1^u\sigma_2^v\sigma_1^w\sigma_2^\epsilon$ and 
$\lambda' = \sigma_1^u\sigma_2^\epsilon\sigma_1^w\sigma_2^v$, and these classes are distinct.  
\item[{\rm (4)}]  If a braid conjugacy class has a representative that is in flype position, then in general it will have at least 3 distinct representatives which are in flype position, related by {\rm (\ref{E:braid representation})}.
\ee
\end{theorem}

{\bf Proof:}  The main theorem in \cite{BM-3-braids} asserts that any link of  braid index $\leq 3$ admits a unique conjugacy class of 3-braid representatives, except in cases (1), (2) and a weak form of (3).  The precise description of the restrictions on $u,v,w$ so that there are precisely two distinct conjugacy classes of 3-braid representatives can be found in Theorem 5 of \cite{KL}. As for {\rm (4)}, see Lemma~\ref{L:conjugate triples} above.$\|$

{\bf Example 2:}  An interesting example is the 3-braid $\lambda = \sigma_1^5\sigma_2^2\sigma_1^3\sigma_2^{-1}$, which admits a negative  flype to $\lambda' = 
\sigma_1^5\sigma_2^{-1}\sigma_1^3\sigma_2^2.$  Rewriting the latter, and using the relations in the braid group $B_3$ to change the representing word, we see that our closed 3-braid $\lambda'$ also admits a positive flype:
$$ \lambda' = (\sigma_1^5)(\sigma_2^{-1}\sigma_1^3\sigma_2)(\sigma_2) = 
(\sigma_1^5)(\sigma_1\sigma_2^3\sigma_1^{-1})(\sigma_2)  = 
\sigma_1^6\sigma_2^3\sigma_1^{-1}\sigma_2,  \ \   \lambda^{\prime\prime} =
\sigma_1^6\sigma_2\sigma_1^{-1}\sigma_2^3,$$
so one might expect there are three conjugacy classes, contradicting the assertion in \cite{BM-3-braids} that there are only two. However it turns out that $\lambda'$ and $\lambda^{\prime\prime}$ are conjugate to one-another, but not conjugate to $\lambda$, so there are only two.

{\bf Example 3:}  By Lemma 1, $\sigma_1^{-5}\sigma_2^{3}\sigma_1^{-3}\sigma_2^{-1}\approx 
\sigma_1^{2}\sigma_2^{-2}\sigma_1^{-5}\sigma_2^{-1} \approx 
\sigma_1^{-3}\sigma_2^{-4}\sigma_1^{2}\sigma_2^{-1}$.  Later we will use the fact that conjugate braids which both admit flypes need not have the same braid crossing number $c_b$. In this case the first has braid crossing number  12 and the other two have crossing number $10$. There are examples where  all three have the same crossing number, and examples where one has minimal crossing number and the other two have higher crossing number.  This will be important to us later, when we enumerate the low-crossing number examples.

Our next result is a  classification theorem for transversal knots that are closed 3-braids, Since we will always be working with negative flypes, we simplify the notation, replacing $ \sigma_1^u\sigma_2^v\sigma_1^w\sigma_2^{-1}$ by $(u,v,w)$:

\begin{theorem}
\label{T:transversal closed 3-braids}  
Let $\mathcal{TK}$ be a transversal knot type which is represented by a closed 3-braid $\kappa$.  
Then its topological knot type is determined uniquely by the conditions given in Theorem \ref{T;classification of links that are closed 3-braids}.  Moreover, $\mathcal{TK}$ is transversally simple except as follows: 
\be
\item [{\rm (1)}]  Choose pairwise distinct integers $(u,v-1,w)$.   
Assume that  $|u|, |v|, |w| \geq 2$ and that $u,w \not= -2$, also $v\not= +2$.  Assume that either (i) $u$ odd and $v,w$ even, or (ii) $w$ odd and $u,v$ even, or (iii) $u,v,w$ all odd.  Then the 3-braids
$\kappa_1 = \sigma_1^u\sigma_2^v\sigma_1^w\sigma_2^{-1}$ and 
$\kappa_2 = \sigma_1^u\sigma_2^{-1}\sigma_1^w\sigma_2^v$     close to distinct transversal knot types $\mathcal{TK}_1, \mathcal{TK}_2$ which have the same topological knot type and the same Bennequin invariant.   That is, they are not transversally simple.
\item [{\rm (2)}]  Call a triple $(u,v,w)$ that satisfies the restrictions given in {\rm (1)} above an {\it admissible triple}.  By Lemma~\ref{L:conjugate triples} we have $(u,v,w)\approx (v-1, w+1,u)\approx (w,u+1,v-1)$.  Then $(u,v,w)$ is admissible if and only if any one of the other two triples is admissible.
\ee
\end{theorem}

{\bf Proof :}  (1) The listed cases are a subset of those in Case 3 of  Theorem \ref{T;classification of links that are closed 3-braids}.  The subset is determined by (a) considering only knots, and (b) considering only 3-braids whose conjugacy class contains a braid that admits a negative flype, and (c) among those, eliminating all classes in which the conjugacy class also admits a positive flype.  The the parity restrictions (i),(ii),(iii) are needed to insure that the permutation  associated to $\kappa$ is a 3-cycle, i.e. $\kappa$ is a knot.  The condition  $\epsilon=-1$ means that  $\kappa$ admits a negative flype.  See Theorem 6 of \cite{KL} for a proof that the conditions $u\not=1$ and $w\not= 1$ and $v\not= 2$ are necessary and sufficient so that $\kappa$ does not also admit a positive flype, for if it did then there would be a transversal isotopy from the transversal knot determined by $\kappa$ to the transversal knot determined by $\kappa'$.  

(2) Replace $u,v,w$ by $u,v-1,w$ for the moment, so that in the new triplet all entries are distinct. Passing to Lemma~\ref{L:conjugate triples}, one now verifies that    the modified triples are cyclic permutations  of $u,v-1,w$, therefore the original triple contains 3 distinct numbers if and only if its conjugate triples do too.  The same is true for the parity conditions.  As for the conditions
$|u|, |v|, |w| \geq 2$ and $u,w \not= -2$, $v\not= +2$, one simply checks the various cases, one at a time. $\|$ 

{\bf Examples 4: An enumeration of all the non-transversally simple knots of braid index 3 which have low crossing number:}  We give a table of all of the examples of non transversally simple knots of braid index 3 which have braid crossing number $c_b = |u|+|v|+|w| + 1\leq 12$.  Note that $c_b$ is not the same as minimum crossing number for all possible projections of the knot associated to our closed braid (although for all but one of the examples in Table 1 the two are actually the same), rather it is crossing number for closed 3-braid representatives which are in flype position. 

To assemble the data in Table 1 we began by listing all possible numbers $(\pm|u|,\pm|v|,\pm|w|)$   with $c_b=|u|+|v|+|w|+1$, deleting those that do not satisfy the  parity restrictions (i),(ii) or (iii) and the other constraints of Theorem~\ref{T:transversal closed 3-braids}.  Next we used Knotscape \cite{Knot} to sort the survivors  into classes which had the same  topological knot type.  For each topological knot type, some of the examples had non-minimal braid crossing number, so we deleted those. For example, $(-3,4,2)$ appeared as an example with $c_b=10$, but (by Lemma~\ref{L:conjugate triples})  its conjugacy class also has an 8-crossing example, $(3,-2,2)$, so we deleted $(-3,4,2)$.  After that, for each topological knot type we divided the representing closed braids into conjugacy classes belonging to the braids before and after the flype.  In every case it turned out that there were two distinct  conjugacy classes in $B_3$, as the main theorem in \cite{BM-3-braids} says there are.  Finally, we chose exactly one representative from each conjugacy class, whenever it happened that the list of survivors had more than one representative.   

The examples are given by specifying,  for each pair of transversal knot types, the topological knot type $\mathcal{K}$ in the standard tables, the  Bennequin number $\beta = u+v+w-4$, the minimum closed braid crossing number $c_b=|u|+|v|+|w|+1$, and the two closed 3-braid representatives ${\mathcal {TK}}_1, {\mathcal {TK}}_2$. The closed braids are given as triplets $(u,v,w)$ and $ (w,v,u)$, i.e. as negative flype-related braids $\sigma_1^u\sigma_2^v\sigma_1^w\sigma_2^{-1}$ and $\sigma_1^w\sigma_2^v\sigma_1^u\sigma_2^{-1}$.   

\

\centerline{
\begin{tabular}{|c | c | c | c |c |  }
\hline
${\mathcal K}$ & $\beta$ & $c_b$ &$\mathcal{TK}_1$ & $\mathcal{TK}_2$ \\
\hline
\hline
$8_{a3}$ & -1 &  8 & (3,-2,2) & (2,-2,3) \\ \hline
$10_{a15}$ &  +1  & 10 &(5,-2,2) & (2,-2,5) \\\hline
$10_{n16}$ &  -9       &  10    & (-5,-2,2) & (2,-2,-5)\\\hline
$10_{n  26}$ & -7 &  10  & (3,-2,-4) & (-4,-2,3) \\\hline
$10_{a 41}$ & +1 &  10  & (3,-2,4) & (4,-2,3)\\\hline
$10_{a79}$ & -3 &   10 & (3,-4,2) & (2,-4,3)\\ \hline
$11_{a240}$ & +7  &   12 & (5,3,3)& (3.3.5)  \\ \hline
$12_{a 146}$ & +3 &   12 & (7,-2,2)& (2,-2,7)  \\\hline
$12_{n 234}$ & -11 &    12 & (-7,-2,2)&  (2,-2,7) \\\hline
$12_{a 369}$ & +3 &   12  & (5,-2,4) &  (4,-2,5) \\\hline
$12_{n  466}$ & -9 &   12  & (-5,-3,3) & (3,-3,-5)  \\\hline
$12_{n  467}$ & -7 &   12  & (-3,-4,4) & (4,-4,-3)  \\\hline
$12_{n  468}$ & -5 &    12  & (-3,-3,5) & (5,-3,-3) \\\hline
$12_{n  472}$ & -15 &    12  & (-3,-4,-4)& (-4,-4,-3)  \\\hline
$12_{n  570}$ & -9 &   12  & (-3,-5,3)& (3,-5,-3) \\\hline
$12_{a 576}$ & +3 &  12  & (3,-2,6) & (6,-2,3)  \\\hline
$12_{a  835}$ &-5 &   12  & (3,-6,2)& (2,-6,3) \\\hline
$12_{a  878}$ & -1 &    12  & (5,-4,2) & (2,-4,5)  \\\hline
$12_{a  1027}$ & -1 &    12  & (3,-4,4) & (4,-4,3)  \\\hline
$12_{a  1233}$ & +1 &   12   &(5,-3,3) & (3,-3,5)  \\
\hline
\end{tabular} }

\

\centerline{{\bf Table 1:  A list of all examples of non-transversally simple }}

\centerline{{\bf  closed 3-braids with braid crossing number $c_b\leq 12$}}

\

{\bf Final Remarks:}  If one is handed a transversal  knot which is represented by a closed 3-braid $\kappa$, the braid $\kappa$ will in general not be in the nice form $\sigma_1^u\sigma_2^v\sigma_1^w\sigma_2^{\epsilon}$.  Of course, a generic example has a unique conjugacy class of 3-braid representatives, but how can we recognize algorithmically whether our braid  is one of the exceptional cases in Theorem \ref{T;classification of links that are closed 3-braids}?  For that purpose one needs to use one of the known solutions to the conjugacy problem in the 3-strand braid group, and to have at hand a list of the classes which admit positive and negative flypes.  Fortunately, Ko and Lee solved that problem, in \cite{KL}, which was written after \cite{BM-3-braids} was published.  Using the solution to the conjugacy problem in the 3-string braid group that was discovered by Peijun Xu in  \cite{Xu}, they define representative symbols which characterize conjugacy class.  The flype-admissible 3-braids are divided into 5 types, but only 3 are needed because if a braid is type $k$, then its inverse is type $6-k$. Proposition 10 of \cite{KL} gives a list of representative  symbols  of flype-admissible 3-braids of types 1,2,3.  To recognize whether a closed 3-braid admits a flype, it suffices to compute the representative symbol, and to check whether it or the symbol of the inverse is in the list given in Proposition 10 of \cite{KL}, also see Theorem 12 of \cite{KL}. 

{\bf Acknowledgment:} We thank Sang-Jin Lee, Lenny Ng, Peter Ozsvath and Dylan Thurston for useful comments.

 \end{document}